\documentclass{article}
\usepackage{amsthm,amsfonts,amssymb,amsmath}
\usepackage{enumerate}
\usepackage[all]{xy}
\usepackage[colorlinks=true]{hyperref}

\newcommand{\BC}{{\mathbb {C}}}

\newcommand{\BR}{{\mathbb {R}}}

\newcommand{\CC}{\mathbb{C}}
\newcommand{\RR}{\mathbb{R}}
\newcommand{\ZZ}{\mathbb{Z}}
\newcommand{\QQ}{\mathbb{Q}}

 \newcommand{\CCC}{{\mathcal {C}}}

\newcommand{\CI}{{\mathcal {I}}}

\newcommand{\CL}{{\mathcal {L}}}
\newcommand{\CM}{{\mathcal {M}}}

\newcommand{\CO}{{\mathcal {O}}}
\newcommand{\CP}{{\mathcal {P}}}
\newcommand{\CQ}{{\mathcal {Q}}}

\newcommand{\CX}{{\mathcal {X}}}
\newcommand{\CY}{{\mathcal {Y}}}

\renewcommand{\div}{{\mathrm{div}}}

\newcommand{\Ind}{{\mathrm{Ind}}}

\newcommand{\Pic}{\mathrm{Pic}}

\newcommand{\vol}{{\mathrm{vol}}}

\newcommand{\sm}{{\mathrm{sm}}}

\renewcommand{\d}{\textnormal{d}}

\renewcommand{\i}{\mathrm{i}}

\newcommand{\wt}{\widetilde}
\newcommand{\wh}{\widehat}

\newcommand{\pair}[1]{\langle {#1} \rangle}

\newcommand{\ds}{\displaystyle}

\newcommand{\ol }{\overline}

\newcommand{\lra}{\longrightarrow}

\newcommand{\CLL}{{\overline{\mathcal L}}}
\newcommand{\CMM}{{\overline{\mathcal M}}}

\newcommand{\CAA}{{\overline{\mathcal A}}}
\newcommand{\CBB}{{\overline{\mathcal B}}}

\newcommand{\CPP}{{\overline{\mathcal P}}}

\newcommand{\OK}{{\overline{K}}}

\newcommand{\lb}{\mathcal{L}}              
\newcommand{\lbb}{\overline{\mathcal{L}}}  
\newcommand{\mbb}{\overline{\mathcal{M}}}  
\newcommand{\abb}{\overline{\mathcal{A}}}  

\newtheorem{thm}{Theorem}[section]

\newtheorem{lem}[thm]{Lemma}

\theoremstyle{remark}
\newtheorem{remark}[thm]{Remark}

\theoremstyle{definition}

\newcommand{\VCV}{\Vert\cdot\Vert}

\newcommand{\Ar}{\textnormal{Ar}}
\newcommand{\achi}{\hat{\chi}}

\newcommand{\gs}[1]{\Gamma(\mathcal{X}, #1)}        
\newcommand{\supnorm}[1]{\|#1\|_{\mathrm{sup}}}    

\DeclareMathOperator{\Spec}{\textnormal{Spec}}
\DeclareMathOperator{\dif}{\textnormal{d}\!}
\DeclareMathOperator{\ord}{\textnormal{ord}}

\begin{document}

\title{On Vojta's proof of the Mordell conjecture}
\author{Xinyi Yuan}
\maketitle
\tableofcontents

\section{Introduction}

The celebrated Mordell conjecture, proved by Faltings, is the following theorem. 

\begin{thm}[Mordell conjecture]
Let $C$ be a smooth, projective, and geometrically connected curve of genus $g>1$ over a number field $K$. 
Then $C(K)$ is finite. 
\end{thm}

The conjecture was raised by Mordell when he proved the Mordell--Weil theorem for elliptic curves over $\QQ$ in 1922.
 Faltings' proof of the conjecture in \cite{Fal83} signified a milestone in the history of Diophantine geometry. 
By an analogue to the Thue--Siegel--Roth theorem in Diophantine approximation, Vojta \cite{Voj91} gave a new proof of the Mordell conjecture in 1991. 
In terms of $p$-adic Hodge theory, Lawrence--Venkatesh \cite{LV20}
gave a third proof of the Mordell conjecture.

Each  of these three proofs has its own valuable feature. For example, Vojta's proof was extended by Faltings \cite{Fal91} to prove the Bombieri--Lang conjecture for subvarieties of abelian varieties.
Moreover, Vojta essentially proved a deep inequality concerning distribution of rational points of large heights. Combining with the uniform Bogomolov conjecture recently proved by Dimitrov--Gao--Habegger \cite{DGH21} and K\"uhne \cite{Kuh}, Vojta's inequality implies a uniform bound on the number of rational points in terms of the genus and the Mordell-Weil rank. 

Vojta's proof depends heavily on high-dimensional Arakelov geometry developed by Arakelov \cite{Ara74}, Deligne \cite{Del}, and Gillet--Soul\'e \cite{GS90, GS92}.
Bombieri \cite{Bom90} (based on \cite{Fal91}) wrote a relatively elementary variant of Vojta's proof, which replaced the use of Arakelov geometry by those of Siegel's lemma and classical height theory, at the cost of involving lots of less important extra objects and estimates. We refer to the textbooks 
\cite{HS00, BG,IKM22} for Bombieri's version of the proof.

The goal of this paper is to introduce a new variant of Vojta's original proof.
The most technical part of Vojta's original proof is the application of the arithmetic Riemann--Roch theorem of Gillet--Soul\'e \cite{GS92} to construct a small section, which also involved an estimate of some analytic torsion. 
Our treatment here replaces the application of the arithmetic Riemann--Roch theorem by a quick application of the arithmetic Siu inequality of Yuan \cite{Yua08}.
So our treatment is still based on Arakelov geometry. 

Note that the arithmetic Siu inequality in \cite{Yua08} is based on Gillet--Soul\'e's arithmetic Hilbert--Samuel formula and Zhang's arithmetic 
Nakai--Moishezon criterion in \cite{Zha92, Zha95}. 
The arithmetic Hilbert--Samuel formula is originally proved by Gillet--Soul\'e \cite{GS92} as a consequence of their arithmetic Riemann--Roch theorem, but it also has a direct proof by Abbes--Bouche \cite{AB}. 
In this sense, our new proof avoids using the arithmetic Riemann--Roch theorem. 

The proofs in \cite{Voj91, Bom90} emphasize the role of the bounds of the index, which is analogous to the proof of the Thue--Siegel--Roth theorem. 
However, in our exposition, we will emphasize more on the role of the bounds of the height and its relation with small sections and base points, which is closer to the spirit of Arakelov geometry. 

Most parts of our treatment are essentially duplicated from Vojta's original proof, but we still include a complete exposition in this paper for convenience of readers.

\subsubsection*{Acknowledgments}
The author would like to thank Zheng Xiao and Chunhui Liu, who taught a course on Roth's theorem and Vojta's inequality in Peking University in the fall of 2022. Motivated by their lectures, the author obtained the variant of Vojta's proof in the current paper. 

The author would also like to thank the support of the China--Russia Mathematics Center. The author is supported by grants NO. 12250004 and NO. 12321001
from the National Science Foundation of China, and by the Xplorer Prize from the New Cornerstone Science Foundation.

\section{Vojta's inequality}

The goal of this section is to introduce Vojta's inequality, Mumford's equality, and  another version of Vojta's inequality using Arakelov geometry. In the end, we sketch the plan of our proof of Vojta's inequality. 

\subsection{Vojta's inequality}

Let $K$ be a number field. Let $C$ be a smooth, projective, and geometrically integral curve of genus $g>1$ over $K$. 
Denote by $J$ the Jacobian variety of $C$ over $K$.  
Denote by $[m]:J\to J$ the morphism given by multiplication by $m\in \ZZ$.

First we recall some basic facts on the N\'eron--Tate heights on $C$ and $J$. Since 
$$ [2g-2]\colon\Pic^0(C_{\ol K})\longrightarrow\Pic^0(C_{\ol K})$$
is surjective, there is a line bundle $\alpha$ on $C_{\ol K}$ such that $(2g-2)\alpha$ is isomorphic to the canonical sheaf $\omega_C$. 
Replacing $K$ by a finite extension if necessary, we can assume that $\alpha$ is actually a line bundle on $C$. 

Consider the Abel--Jacobi embedding
$$
i_\alpha: C\lra J, \quad  x\longmapsto (x)-\alpha.
$$
Recall that the theta divisor on $J$ is given by
$$
\theta_\alpha= \underbrace{i_\alpha(C)+\dots + i_\alpha(C)}_{g-1 \text{ copies}}.$$
It is well-known that $\theta_\alpha$ is ample and gives a principal polarization of $J$. 
By \cite[p. 74, eq. (1)]{Ser89}, $\theta_\alpha$ is symmetric in the sense that $[-1]\theta_\alpha$ is linearly equivalent to $\theta_\alpha$. 

Denote $J(\ol K)_\RR=J(\ol K)\otimes_\ZZ \RR$. 
It is well-known that the N\'eron--Tate height 
$$ \hat h_{\theta_\alpha}\colon J(\ol K)\lra \RR$$
is quadratic and induces a positive definite quadratic form
$$ \hat h_{\theta_\alpha}\colon J(\ol K)_\RR\lra \RR.$$
 See \cite[\S3.3, \S3.8]{Ser89} for example, and we normalize Weil heights over number fields as in \cite[\S2.2]{Ser89}.
Denote the N\'eron--Tate height pairing
$$\pair{P,Q}\colon J(\overline{K})\times J(\overline{K})\longrightarrow \BR$$
by 
$$ \pair{P,Q}=\frac12 \big(\hat{h}_{\theta_\alpha}(P+Q)-\hat{h}_{\theta_\alpha}(P)-\hat{h}_{\theta_\alpha}(Q)\big), \quad P, Q\in J(\ol K).$$
To do Euclidean geometry, we define 
$$|P| := \sqrt{\hat{h}_{\theta_\alpha}(P)}, \quad P\in J(\ol K).$$
We have $|P| = \sqrt{\pair{P,P}}$ by definition. 
These definitions extend to $P, Q\in J(\ol K)_\RR$ naturally. 
We will apply $ |\cdot|$ and $\pair{\cdot,\cdot}$ to $C(\ol K)$ via the embedding $i_\alpha: C\to J$. 

Now we are ready to introduce Vojta's inequality, which is essentially due to Vojta \cite{Voj91}. 

\begin{thm}[Vojta's inequality] \label{vojta}
For any real number $\ds\mu>\frac{1}{\sqrt g}$, there exist positive constants $A_1, A_2$ such that for all points $P_1,P_2\in C(\ol K)$ satisfying 
$$|P_1|>A_1, \quad  
|P_2|>A_2|P_1|, 
$$
we have
$$
\pair{P_1,P_2}< \mu\cdot |P_1|\cdot |P_2|.
$$
\end{thm}

Combined with the Mordell--Weil theorem for the Jacobian variety, Vojta's inequality implies the Mordell conjecture without much effort. 
In fact, the finite-dimensional real vector space $V=J(K)_\RR$, endowed with the inner product $\pair{\cdot,\cdot}$, is isometric to a Euclidean space $\RR^r$.
For any two vectors $x,y\in V$, denote by $\angle(x,y)$ the angle between these two vectors, which is set to be 0 if $x=0$ or $y=0$.  
Take a real number $\mu$ in the open interval $(1/\sqrt g, 1)$, and denote $\theta=\arccos \mu$.  
The conclusion
$$
\pair{P_1,P_2}< \mu\cdot |P_1|\cdot |P_2|
$$
implies  
$$ \cos \angle(P_1,P_2)< \mu, \quad \angle(P_1,P_2)> \theta.$$

For every nonzero vector $x\in V$, consider the cone
$$
\mathrm{cone}(x):=\{y\in V: \angle(x,y)\leq \theta/2\}.
$$
Then $V$ is covered by finitely many cones 
$\mathrm{cone}(x_1), \dots, \mathrm{cone}(x_n)$, as a consequence of the compactness of the unit sphere of $V$.
For every $i$, denote by $\Sigma_i$ the set of points of $C(K)$ whose images in $V$ lie in $\mathrm{cone}(x_i)$. 
It suffices to prove that every $\Sigma_i$ is finite.
Note that any two points of $\Sigma_i$ have an angle at most $\theta$.
Apply Theorem \ref{vojta} to this situation. 
If every $P\in \Sigma_i$ satisfies $|P|\leq A_1$, then $\Sigma_i$ is finite by the Northcott property.
If there exists $P\in \Sigma_i$ with $|P|> A_1$, then Vojta's inequality implies 
$|P'|\leq A_2|P|$ for any 
$P'\in \Sigma_i$, and thus $\Sigma_i$ is still finite by the Northcott property.

\subsection{Mumford's equality}

To prove Vojta's inequality, we first apply Mumford's equality to convert the inequality to a lower bound of Weil heights of some relevant points. 

Let $K$ be a number field. Let $C$ be a smooth projective and geometrically connected curve of genus $g>1$ over $K$. 
As above, let $\alpha$ be a line bundle on $C$ with 
$(2g-2)\alpha\cong \omega_{C/K}$. 

Let $X=C^2$, and let $\Delta$ be the diagonal of $C^2$. 
Denote by $p_1,p_2:X\to C$ the two projections. Define the \emph{Vojta line bundle} to be
$$
L_0:= (\delta_1-1) p_1^* \alpha+(\delta_2-1) p_2^* \alpha+ \CO(\Delta), 
$$
where $\delta_1,\delta_2$ are positive rational numbers. The Vojta line bundle $L_0$ is a $\QQ$-line bundle on $X$, i.e. an element of $\Pic(X)_\QQ=\Pic(X)\otimes_\ZZ\QQ$.
Here we write tensor products of line bundles additively. 

We will see that it is often more convenient to write
$$
L_0= \delta_1\, p_1^* \alpha+\delta_2\, p_2^* \alpha+ \CO(\Delta)^\circ,
$$
where
$$
\CO(\Delta)^\circ:=\CO(\Delta)-p_1^* \alpha-p_2^* \alpha.
$$
is the  the modified diagonal. For example, it gives convenience in computing intersection numbers by
$$
\CO(\Delta)^\circ \cdot p_1^* \alpha
=\CO(\Delta)^\circ \cdot p_2^* \alpha
=0.
$$

Let $h_{p_1^* \alpha}, h_{p_2^* \alpha}, h_\Delta$ be Weil heights on $X(\ol K)$ with respect to the line bundles $p_1^*\alpha, p_2^* \alpha, \CO(\Delta)$.
Take a Weil height function
$$
h_{L_0}: X(\OK)\lra \RR
$$
with respect to the line bundle $L_0$ on $X$ by
$$
h_{L_0}= (\delta_1-1) h_{p_1^* \alpha}+(\delta_2-1) h_{p_2^* \alpha}+ h_\Delta. 
$$
Note that $h_{p_1^* \alpha}$, $h_{p_2^* \alpha}$, and $h_\Delta$ are only determined up to constants, but we choose them once for all. 

\begin{thm}[Mumford's equality] \label{mumford}
For any $P_1,P_2\in C(\ol K)$, 
$$
h_{L_0}(P_1,P_2)
=\frac{1}{g}\delta_1  |P_1|^2
+\frac{1}{g} \delta_2 |P_2|^2 
-2 \pair{P_1,P_2}+O(\delta_1+\delta_2+1).
$$
Here the coefficient in the error term has a bound independent of $(P_1,P_2, \delta_1, \delta_2)$.
\end{thm}

\begin{proof}
Let us start with the pull-back formulas. 
Denote by 
$$j_\alpha=(i_\alpha, i_\alpha):C\times C\lra J\times J$$ 
the natural morphism. 
Denote by $p_1', p_2':J \times J\to J$ the projections.
Denote by $P$ the Poincar\'{e} line bundle on $J\times J$. 
Then the following are true:
\begin{enumerate}[(1)]
\item $i_\alpha^*\CO(\theta_\alpha)=g\alpha$ in $\Pic(C)$. 
\item $j_\alpha ^*P = \CO(\Delta)^\circ$ 
in $\Pic(C\times C)$. 
\end{enumerate}
The first identity is \cite[p. 75, eq. (2)]{Ser89}, the second identity is \cite[p. 76, eq. (4)]{Ser89}. 
All these require the condition $(2g-2)\alpha=\omega_{C/K}$. 

Apply the functoriality of Weil heights for pull-backs of line bundles as in \cite[p. 23]{Ser89} to the above identities. We have for $k=1,2$, 
$$
h_{p_k^*\alpha}(P_1,P_2)=h_{\alpha}(P_k)+O(1)=
\frac1g h_{\theta_\alpha}(i_\alpha(P_k))+O(1)
=\frac1g |P_k|^2+O(1),
$$
and 
$$
h_{ \CO(\Delta)^\circ}(P_1,P_2)
=h_{P}(i_\alpha(P_1), i_\alpha(P_2))+O(1)
=2\pair{P_1, P_2}+O(1).
$$ 
Here the last inequality follows from \cite[p. 39, Theorem]{Ser89},
and the factor 2 comes from different normalization of the pairing. 
\end{proof}

\begin{remark}
We can have a precise equality of the above theorem without the error $O(1)$
in terms of admissible pairings of Zhang \cite{Zha93} (or the reformulation of \cite[Appendix A]{Yua21} in terms of adelic line bundles of Zhang \cite{Zha95b}).
\end{remark}

In the case $\delta_1=\delta_2=1$, Mumford's equality gives
$$
h_{\Delta}(P_1,P_2)
=\frac{1}{g}  |P_1|^2
+\frac{1}{g}  |P_2|^2 
-2 \pair{P_1,P_2}+O(1).
$$
As $\Delta$ is effective, its Weil height function satisfies
$$
h_{\Delta}(P_1,P_2)\geq O(1), \quad P_1\neq P_2. 
$$
This implies an inequality  
$$
\frac{1}{g}  |P_1|^2
+\frac{1}{g}  |P_2|^2 
-2 \pair{P_1,P_2}\geq O(1), \quad P_1\neq P_2. 
$$
This inequality is due to Mumford  \cite[\S3, Cor. 1]{Mum65}. It 
gives a non-trivial bound on the angle between $P_1$ and $P_2$ if 
$|P_1|$ and $|P_2|$ are large with a quotient close to 1. See also \cite[\S5.7]{Ser89}.

To illustrate Vojta's idea to use Mumford's equality to prove Vojta's inequality, let us assume some situations which are too good to be true but might eventually be modified to be true. 
Let $\delta_1, \delta_2$ be positive rational numbers.
A basic application of the Riemann--Roch theorem on $X$ proves that Vojta's line bundle 
$$
L_0= \delta_1 p_1^* \alpha+\delta_2 p_2^* \alpha+ \CO(\Delta)^\circ
$$
has volume
$$
\vol(L_0)\geq L_0^2=2(\delta_1\delta_2-g).
$$
In particular, $L_0$ is big if $\delta_1\delta_2>g$. See Theorem \ref{geometric volume} for this result.
We will assume that $\delta_1\delta_2>g$ throughout this paper.
It follows that $h_{L_0}(P)$ is bounded below for $P=(P_1,P_2)\in X(\ol K)$ lying outside of the base locus of positive multiples of $L_0$. See \cite[\S2.10, Theorem]{Ser89}. 
By Mumford's equality, this gives the relation 
$$
\frac{1}{g}\delta_1  |P_1|^2
+\frac{1}{g} \delta_2 |P_2|^2 
-2 \pair{P_1,P_2}
\geq 
O(\delta_1+\delta_2+1)
$$
for $P\in X(\ol K)$ outside the base locus.

For the sake of demonstrating ideas, we ignore the base locus and uniformity of the error, so we simply assume that under the condition $\delta_1\delta_2>g$, we have
$$
\frac{1}{g}\delta_1  |P_1|^2
+\frac{1}{g} \delta_2 |P_2|^2 
-2 \pair{P_1,P_2}
\geq 
0
$$
for all $P_1,P_2\in C(\ol K)$.
For fixed $(P_1,P_2)$, we want to choose $\delta_1, \delta_2$ so that the inequality is as optimal as possible. 
We also fix $\lambda=\delta_1\delta_2-g>0$. 
The most optimal case amounts to taking the minimal value of the left-hand side of the inequality under the constraint $\delta_1\delta_2=g+\lambda$. 
By the inequality between the arithmetic mean and the geometric mean, the left-hand side of the inequality is minimal if 
$$
\delta_1=  \sqrt{g+\lambda} \frac{|P_2|}{|P_1|}, \qquad
\delta_2=\sqrt{g+\lambda} \frac{|P_1|}{|P_2|}.
$$
In this case, the inequality becomes
$$
\frac{\sqrt{g+\lambda}}{g} |P_1|\cdot  |P_2|
\geq  \pair{P_1,P_2}.
$$
This would imply Vojta's inequality.

However, this ideal situation is not practical by many obstacles.
First, the point $(P_1,P_2)$ might lie in the base locus of the line bundle $L_0$.
Second, we need uniformity of the error term when varying $(P_1,P_2)$ and $(\delta_1, \delta_2)$. 
In the following, we write a precise statement of the inequality, and then introduce more ideas to overcome these two issues. 

\subsection{A precise statement}
Let $\delta_1,\delta_2$ be positive rational numbers. 
As above, denote 
$$
L_0= (\delta_1-1) p_1^* \alpha+(\delta_2-1) p_2^* \alpha+ \CO(\Delta). 
$$
This gives a Weil height 
$$
h_{L_0}:= (\delta_1-1) h_{p_1^* \alpha}+(\delta_2-1) h_{p_2^* \alpha}+ h_\Delta 
$$
on $X=C^2$ as above. 
The estimate we can achieve is as follows. 

\begin{thm}\label{vojta2}
Let $\delta_1,\delta_2$ be positive rational numbers satisfying
\begin{equation*} 
\delta_1>2g\delta_2,  \quad  g<\delta_1\delta_2<g+\frac14.
\end{equation*}
Let $P_1,P_2\in C(\ol K)$ be points.
Then
$$
h_{L_0}(P_1,P_2)
\geq 
-2\, c(\delta_1,\delta_2)\cdot (\delta_1|P_1|^2+\delta_2|P_2|^2)
+O(\delta_1/(\delta_1\delta_2-g)).
$$ 
Here the coefficient in the error term is independent of $(P_1,P_2, \delta_1, \delta_2)$, and 
$$
c(\delta_1,\delta_2)=\sqrt{\frac{2}{g}(\delta_1\delta_2-g)+(2g-1)\frac{\delta_2}{\delta_1}}.
$$
\end{thm}

By Mumford's equality, we can prove that Theorem \ref{vojta2} implies Theorem \ref{vojta}. 
In fact, by Theorem \ref{mumford}, Theorem \ref{vojta2} gives
$$
\left(\frac{1}{g}
+2\, c(\delta_1,\delta_2)\right)  (\delta_1|P_1|^2+\delta_2|P_2|^2)
-2 \pair{P_1,P_2}
+ \frac{c_0\delta_1}{\delta_1\delta_2-g} \geq 0.
$$
Here $c_0$ is a  non-negative constant independent of $(P_1,P_2, \delta_1,\delta_2)$.
By approximation, the inequality also holds for real numbers $\delta_1,\delta_2$ satisfying the condition of Theorem \ref{vojta2}.

Let $\lambda$ be a real number with $0<\lambda<1/4$. Set
$$
\delta_1=  \sqrt{g+\lambda} \frac{|P_2|}{|P_1|}, \qquad
\delta_2=\sqrt{g+\lambda} \frac{|P_1|}{|P_2|}.
$$
Then $\delta_1\delta_2-g=\lambda>0$.
The requirement $\delta_1>2g\delta_2$ transfers to $|P_2|>\sqrt{2g} |P_1|$, which we assume. We further assume $|P_1|>0$. 
The inequality becomes
$$
2\sqrt{g+\lambda} \left(\frac{1}{g} 
+ 2c(\delta_1,\delta_2)\right) |P_1|\cdot |P_2| 
+ \sqrt{g+\lambda} \frac{c_0|P_2|}{\lambda |P_1|} \geq 2 \pair{P_1,P_2},
$$
or equivalently
$$
\frac{\pair{P_1,P_2}}{|P_1|\cdot |P_2| }\leq 
\sqrt{g+\lambda} \left(\frac{1}{g} 
+ 2c(\delta_1,\delta_2)+ \frac{c_0}{2\lambda |P_1|^2} \right).
$$
Here 
$$
c(\delta_1,\delta_2)=\sqrt{\frac{2}{g}\lambda+(2g-1)\frac{|P_1|^2}{|P_2|^2} }.
$$
For any fixed $\lambda>0$, if $|P_1|^2>\lambda^{-2} c_0$ and 
$|P_2|^2>\lambda^{-1} |P_1|^2$, then the inequality implies
$$
\frac{\pair{P_1,P_2}}{|P_1|\cdot |P_2| }\leq 
\sqrt{g+\lambda} \left(\frac{1}{g} 
+ 2 \sqrt{\frac{2}{g}\lambda+(2g-1)\lambda }+ \frac{\lambda}{2} \right).
$$
As $\lambda\to 0$, the right-hand side converges to $1/\sqrt{g}$. It follows that the right-hand side can be taken to be smaller than any given $\mu>1/\sqrt{g}$.  
This proves Theorem \ref{vojta} by Theorem \ref{vojta2}.

\subsection{The Arakelov-geometric setting}

The Weil heights in Theorem \ref{vojta2} are only well-defined up to bounded functions, which is inconvenient for estimation when varying the line bundles. 
So we are going to realize them as height functions associated to hermitian line bundles.

\subsubsection{Basic Arakelov geometry}

Let us start with some terminology of Arakelov geometry of Arakelov \cite{Ara74} and Gillet--Soul\'e \cite{GS90}. 

Let $\CX$ be an arithmetic variety, i.e. a projective and flat integral scheme over $\Spec \ZZ$. 
Assume that the generic fiber $\CX_\QQ=\CX\times_\ZZ\QQ$ of $\CX$ is smooth. 
Let $\CLL=(\CL,\VCV)$ be a \emph{hermitian line bundle} on $\CX$, i.e., $\CL$ is a line bundle on $\CX$, and $\VCV$ is a (smooth) hermitian metric of $\CL(\CC)$ on the complex manifold $\CX(\CC)$ invariant under the complex conjugation.

For any section $s\in \Gamma(\CX, \lb)_\RR =\Gamma(\CX, \lb)\otimes_{\ZZ}\RR \subset \Gamma(\CX_\CC, \lb_\CC),$ we have the \emph{supremum norm} 
$$\supnorm{s}=\sup_{z\in \CX(\mathbb{C})} \|s(z)\|.$$
 Picking any Haar measure on $\Gamma(\CX, \lb)_\RR$, we have the 
 \emph{arithmetic Euler characteristic}
$$\hat\chi(\lbb)=\log\frac{\vol(B(\lbb))}{\vol(\gs{\lb}_{\RR}/\gs{\lb})},$$
where 
$$B(\lbb)=\{ s\in\gs{\lb}_{\RR}: \supnorm{s}\leq 1\}$$ 
is the corresponding unit ball. It is easy to see that this definition is independent of the choice of the Haar measure.

A hermitian line bundle $\lbb$ on $\CX$ is called \emph{nef} if the following conditions are satisfied:
\begin{enumerate}
\item $\lbb$ is \textit{relatively semipositive}; i.e., the curvature of $\lbb$ is semipositive and $\deg(\lb|_\CY) \geq 0$ for any closed integral curve $\CY$ on any special fibre of $\CX$ over $\mathrm{Spec}(\mathbb{Z})$;
\item $\lbb$ is \textit{horizontally semipositive}; i.e., the arithmetic degree
 $\wh\deg(\lbb|_\CY)\geq 0$ for any one-dimensional horizontal closed integral subscheme $\CY$ of $\CX$.
\end{enumerate}
Note that the second condition in (2) means that $\lb$ is nef on any
special fiber of $\CX$ over $\ZZ$ in the classical sense. 
By an arithmetic version of Kleiman's theorem (cf. \cite[Cor. A.4.3]{YZ}), if $\lbb$ is nef, then $\lbb^{\dim \CY}\cdot \CY\geq 0$ for any closed integral subscheme $\CY$ of $\CX$.


As a convention of this paper, we write tensor products of (hermitian) line bundles additively, so $n(\CLL-\CMM)$ means $(\CLL\otimes \CMM^\vee)^{\otimes n}$. 

Now we have the following arithmetic Siu inequality of Yuan \cite{Yua08}, whose proof is based on the arithmetic Hilbert-Samuel formula of Gillet--Soul\'e \cite{GS92} and the arithmetic Nakai--Moishezon criterion of Zhang \cite{Zha95}.

\begin{thm}\label{bigness}
Let $\lbb$ and $\mbb$ be nef hermitian line bundles on an arithmetic variety $\CX$ of dimension $r$.  Then
$$\achi(n(\lbb-\mbb))
\geq \frac{n^r}{r!}(\lbb^r-r\, \lbb^{r-1}\cdot \mbb)+o(n^r).$$
\end{thm}

\begin{proof}
If $\lbb$ and $\mbb$ are ample, the theorem follows from \cite[Thm. 2.2]{Yua08}.
If $\lbb$ and $\mbb$ are nef, take an ample hermitian line bundle $\abb$ on $\CX$. Apply the theorem to $(m\lbb+\abb, m\mbb+\abb, i(\lbb-\mbb))$ for positive integers $m$ and $i=0,\dots, m-1$. Set $m\to \infty$.
\end{proof}

\subsubsection{Arithmetic models}

We fix the following arithmetic datum:
\begin{enumerate}[(1)]
\item
Let $\CCC$ be a regular (projective and flat) integral model of $C$ over $O_K$. 
\item
Let $\CX$ be a (projective and flat) integral model of $X$ over $O_K$ such that the identification $\CX_K\to C^2$
extends to a morphism $\psi:\CX\to \CCC\times_{O_K}\CCC$ which induces an isomorphism from $\psi^{-1}(\CCC^{\rm sm}\times_{O_K}\CCC^{\rm sm})$ to $\CCC^{\rm sm}\times_{O_K}\CCC^{\rm sm}$.
Here $\CCC^{\rm sm}$ denotes the smooth locus of 
$\CCC$ over $O_K$. 
\item
Denote by $\tilde p_1, \tilde p_2:\CX\to \CCC$ the morphism extending the projections 
$p_1,p_2:X\to C$, which are also the compositions of  $\psi:\CX\to \CCC\times_{O_K}\CCC$ with the projections $p_1,p_2: \CCC\times_{O_K}\CCC\to \CCC$.
\item
Let $\overline\alpha=(\tilde \alpha, \|\cdot\|)$ be a \emph{nef} hermitian $\QQ$-line bundle on $\CCC$ extending $\alpha$.
\item
Let $\overline \Delta$ be an arithmetic divisor 
on $\CX$ extending $\Delta$. 
\end{enumerate}

There are various methods to construct $(\CCC,\CX, \ol\alpha,\ol\Delta)$, since our requirements are very weak.
For example, assume that $C$ has semistable reduction over $O_K$, which can be achieved by a finite base change. Then let $\CCC$ be the minimal regular model of $C$ over $O_K$, and let 
$\ol\omega_{\CCC/O_K}$ be the relative dualizing sheaf with Arakelov's metric at archimedean places. 
By \cite[\S5, Thm. 5]{Fal84}, $\ol\omega_{\CCC/O_K}$ is nef. 
Note that this result is not essential for our purpose, since we could always change the hermitian metric by a small constant multiple to make $\ol\omega_{\CCC/O_K}$ nef. 
Take $\ol\alpha$ to be the unique extension of $\alpha$ such that
$(2g-2)\ol\alpha=\ol\omega_{\CCC/O_K}$.
Then we set $\CX$ to be the blowing-up of $\CCC\times_{O_K} \CCC$ along its 
(0-dimensional) non-smooth locus, which is known to be regular by an easy explicit computation (cf. \cite[\S3.1]{Zha10}).
Set $\ol\Delta$ to be the Zariski closure of $\Delta$ in $\CX$ with an arbitrary Green function.

Finally, we define Vojta's  \emph{hermitian $\QQ$-line bundle} $\CLL_0$ on $\CX$ by 
$$
\CLL_0:= (\delta_1-1) \tilde p_1^* \ol\alpha+(\delta_2-1) \tilde p_2^* \ol\alpha+ \CO(\ol\Delta). 
$$
Recall the \emph{height function} 
$$
h_{\CLL_0}: X(\OK) \lra \RR
$$
by 
$$
h_{\CLL_0}(x)=\frac{1}{\deg(x)} \wh\deg(\CLL_0|_{\bar x}),
$$
where $\bar x$ is the multi-section of $\CX$ over $\Spec O_K$ corresponding to $x$, and $\deg(x)$ is the degree over $\bar x$ over $\Spec O_K$.  
For convenience, we normalize the height function by 
$$
h_{\CLL_0}(x)_{\QQ}=\frac{1}{[K:\QQ]}h_{\CLL_0}(x),\quad x\in X(\OK).
$$
Then the following theorem is equivalent to Theorem \ref{vojta2}. 

\begin{thm}\label{vojta3}
Let $\delta_1,\delta_2$ be positive rational numbers satisfying
\begin{equation*} 
\delta_1>2g\delta_2,  \quad  g<\delta_1\delta_2<g+\frac14.
\end{equation*}
Let $P_1,P_2\in C(\ol K)$ be points.
Then
$$
h_{\CLL_0}(P_1,P_2)_\QQ
\geq 
-2\, c(\delta_1,\delta_2)\cdot (\delta_1|P_1|^2+\delta_2|P_2|^2)
+O(\delta_1/(\delta_1\delta_2-g)).
$$ 
Here the coefficient in the error term is independent of $(P_1,P_2, \delta_1, \delta_2)$, and 
$$
c(\delta_1,\delta_2)=\sqrt{\frac{2}{g}(\delta_1\delta_2-g)+(2g-1)\frac{\delta_2}{\delta_1}}.
$$
\end{thm}
 
Note that we can assume that $P_1,P_2\in C(K)$ by taking finite base changes, and track that the error term in our proof is uniform in this sense. 
Hence, from now on, we always assume $P_1,P_2\in C(K)$, and omit the part of checking the uniformity under base change.

\subsection{Plan of proof}

Recall that the Mordell conjecture is reduced to Theorem \ref{vojta3}. 
The proof of the theorem will take up the rest of this paper.
Let us sketch the idea here by three major steps. 

The first step is to construct a small section $s$ of $\CLL=d\CLL_0$, i.e. a global section $s$ of $d\CL_0$ with small supremum norm $\|s\|_{\sup}$ over $X(\CC)$. 
Here $d$ is any sufficiently large integer such that $d\CLL_0$ is an integral hermitian line bundle. 
To achieve this, it suffices to have a suitable lower bound of $\achi(d\CLL_0)$.
We apply Yuan's arithmetic Siu inequality (cf. Theorem \ref{bigness}) to achieve this, while Vojta applied the arithmetic Riemann--Roch theorem. 

The second step is to obtain a suitable lower bound of the height $h_{\CLL}(P)$ in terms of the small section $s$. This is easy and standard if $s$ does not vanish at $P$. 
However, it is not practical to expect the non-vanishing, and thus we have to treat the case that $s$ vanishes at $P$. 
The index of $s$ at $P$ (with a proper weight) is a geometric term for ``the vanishing order'' of $s$ at $P$. 
Then we give a lower bound of $h_{\CLL}(P)$ in terms of $-\log\|s\|_{\sup}$ and a contribution from the index. 
The smaller the index is, the better the lower bound is. 
The proof here is mostly to unravel definitions in Arakelov geometry. 

The third step is to obtain an upper bound of the index of $s$ at $P$, which will be used to simplify the lower bound of the height $h_{\CLL}(P)$ to prove Theorem \ref{vojta3}. In history, there are two approaches to such an upper bound. 
Vojta used Dyson's lemma, which is geometric in nature; \cite{Bom90} used Roth lemma, which is arithmetic in nature.
We follow Vojta's approach by Dyson's lemma. 
The key idea to prove Dyson's lemma is that by blowing-up $X$ at $P$,  the index of $s$ at $P$ (with the trivial weight) becomes the multiplicity of the exceptional divisor
in the pull-back of $\div(s)$, and then it can be bounded by intersection numbers of the line bundles involved.

\subsection*{Convention for error terms}

Throughout this paper, our estimates are sensitive to $(P_1,P_2, \delta_1, \delta_2, d)$. As a convention, by an inequality of the form
$$
F\geq G+O(H)
$$
for real number $F,G,H$ with $H>0$, 
we mean 
$$
F\geq G+c H
$$
for a (possibly negative) constant $c$ independent of $(P_1,P_2, \delta_1, \delta_2, d)$. 
To emphasize this situation, we may also say that
 the coefficient in the error term is independent of $(P_1,P_2, \delta_1, \delta_2,d)$.
The convention for $F\leq G+O(H)$ and $F= G+O(H)$ are similar.

\section{Step 1: existence of small section}

Recall that 
 $$
\CLL_0:= (\delta_1-1) \tilde p_1^* \ol\alpha+(\delta_2-1) \tilde p_2^* \ol\alpha+ \CO(\ol\Delta). 
$$
Here $(\CCC,\ol\alpha, \CX,\ol\Delta)$ is the arithmetic model of 
$(C, \alpha, X, \Delta)$ chosen above. 
In particular, $\ol\alpha$ is nef on $\CCC$ by the choice. 
The goal of this step is to prove the following theorem, which is a variant of \cite[Thm. 2.2]{Voj91}.

\begin{thm}[small section]\label{small section}
Let $\delta_1,\delta_2$ be positive rational numbers satisfying
$\delta_1>\delta_2$ and $\delta_1\delta_2>g.$
Then there exists a positive integers $d_0$ depending on $\CLL_0$, such that for any integer $d\geq d_0$ making $d\CLL_0$ an (integral) hermitian line bundle on $\CX$,
there exists a nonzero element
$s\in H^0(\CX, d\CL_0)$ satisfying 
$$
-\log \|s\|_{\sup}
\geq O\left(\frac{d\delta_1^2\delta_2}{\delta_1\delta_2-g}\right).
$$
\end{thm}

\subsection{Hilbert--Samuel type of estimates}

Denote 
$$\alpha_1=p_1^* \alpha, \quad  
\alpha_2=p_2^* \alpha, \quad  
\ol\alpha_1=\wt p_1^* \ol\alpha, \quad 
\ol\alpha_2=\wt p_2^* \ol\alpha.
$$
Recall  
$$
L_0= \delta_1\,  \alpha_1+\delta_2\, \alpha_2+  ( \CO(\Delta)-\alpha_1- \alpha_2)
$$
and
$$
\CLL_0= \delta_1\,  \ol\alpha_1+\delta_2\, \ol\alpha_2
+ ( \CO(\ol\Delta)- \ol\alpha_1- \ol\alpha_2). 
$$

We first have the following geometric estimates, which implies that $L_0$ is big on $X$. Note that we cannot count sections of the $\QQ$-line bundle $L_0$ directly, so take multiples of $L_0$ which are (integral) line bundles on $X$.

\begin{thm} \label{geometric volume}
For positive integers $d$ such that $dL_0$ is a line bundle on $X$, we have 
$$
h^0(dL_0)\geq (\delta_1\delta_2-g)d^2+O(d).
$$
\end{thm}

\begin{proof}
By the Riemann--Roch theorem, for any positive integer $d$ such that $dL_0$ is an integral line bundle on $X$, 
$$
h^0(dL_0)-h^1(dL_0)+h^2(dL_0)= \frac12 d^2 L_0^2-\frac12 d\, L_0\cdot \omega_X+\chi(\CO_X).
$$
By the Serre duality, $h^2(dL_0)=h^0(\omega_X-dL_0)$. 
The line bundle $M=\alpha_1+\alpha_2$ is ample on $X$.
An easy calculation gives
$L_0\cdot M=\delta_1+\delta_2$. 
It follows that for sufficiently large $d$,
the intersection number  $M\cdot(\omega_X-dL_0)<0$, and thus 
$h^0(\omega_X-dL_0)=0$. 
It follows that 
$$
h^0(dL_0)\geq  \frac12 d^2 L_0^2-\frac12 d\, L_0\cdot \omega_X+\chi(\CO_X).
$$
\end{proof}

\begin{remark}
One can also prove that $h^0(dL_0)\leq (\delta_1\delta_2)d^2+O(d)$, but we do not need this fact here.
\end{remark}

The following key estimate can be viewed as an arithmetic counterpart of 
Theorem \ref{geometric volume}. 
We will prove it by Yuan's arithmetic Siu inequality instead of Gillet--Soul\'e's arithmetic Riemann--Roch theorem.

\begin{thm} \label{vol chi}
There is a positive constant $A_{\delta_1,\delta_2}$ depending on $\delta_1, \delta_2$ such that for positive integers $d>A_{\delta_1,\delta_2}$ making $d\CLL_0$ is a hermitian line bundle on $\CX$, we have 
$$\achi(d\CLL_0)\geq \frac{1}{6}d^3(\CLL_0^3)
+O(\delta_1d^3)= O(\delta_1^2\delta_2d^3).$$
\end{thm}

Let us first see how Theorem \ref{geometric volume} and Theorem \ref{vol chi} imply Theorem \ref{small section}
by applying Minkowski's  theorem.
In fact, denote $\epsilon=(\delta_1\delta_2-g)/2$.
By the theorems, 
 for sufficiently large and sufficiently divisible $d$ (relatively to $(\delta_1,\delta_2)$), we have
$$
h^0(dL_0)\geq 
\frac12 d^2 \vol(L_0)-\epsilon d^2
\geq\frac12 d^2(\delta_1\delta_2-g),
$$
and
$$
\achi(d\CLL_0)\geq 
O(\delta_1^2\delta_2d^3).
$$
By Minkowski's theorem, there exists a nonzero section $s\in \Gamma(\CX, d\CL_0)$ such that 
\begin{multline*}
-\log \|s\|_{\sup} 
\geq \frac{1}{[K:\QQ]} \frac{\achi(d\CLL_0)}{h^0(dL_0)}-\log 2\\
\geq \frac{O(d^3\delta_1^2\delta_2)}{d^2(\delta_1\delta_2-g)}-\log 2
=O\left(\frac{d\delta_1^2\delta_2}{\delta_1\delta_2-g}\right).
\end{multline*}
Here the second inequality uses the lower bound of $h^0(dL_0)$ (instead of an upper bound) due to the extreme possibility that $O(d^3\delta_1^2\delta_2)$ might represent a negative number.
This proves Theorem \ref{small section}.

\subsection{Estimate of the arithmetic Euler characteristic}

Now we prove Theorem \ref{vol chi}. 
Recall
$$
\CLL_0=\delta_1\,  \ol\alpha_1+\delta_2\, \ol\alpha_2
+( \CO(\ol\Delta)- \ol\alpha_1- \ol\alpha_2). 
$$
Write 
$$
\CO(\ol\Delta)- \ol\alpha_1- \ol\alpha_2=\CAA-\CBB
$$
for nef hermitian line bundles $\CAA, \CBB$ on $\CX$. 
It follows that
$$
\CLL_0= \CAA'-\CBB 
$$
for 
$$
\CAA'= \delta_1\,  \ol\alpha_1+\delta_2\, \ol\alpha_2+ \CAA.
$$
By the construction, $\ol\alpha$ is nef on $\CCC$, so $\CAA'$ is also nef on $\CX$.
Then we are ready to apply Theorem \ref{bigness}. It gives
$$
\achi(d\CLL_0)\geq \frac{d^3}{6}(\CAA'^3-3\,\CAA'^2\cdot\CBB)+o_{\delta_1,\delta_2}(d^3)
=\frac{d^3}{6}(\CLL_0^3-3\,\CAA'\cdot\CBB^2+\CBB^3)+o_{\delta_1,\delta_2}(d^3).
$$
Here the coefficient of the error term $o_{\delta_1,\delta_2}(d^3)$ depends on $\delta_1,\delta_2$.

By the expressions defining $\CLL_0$, $\CAA'$ and $\CBB$, we see that 
the intersection numbers on the right-hand side
are polynomials of $\delta_1,\delta_2$ of degree at most 3, whose coefficients are given by intersection numbers of $\ol\alpha_1, \ol\alpha_2, \ol\Delta, \CAA, \CBB$ on $\CX$. 
The polynomials appearing in $\CAA\cdot\CBB^2$ and $\CBB^3$ have degrees at most $1$. 
As $\delta_1>\delta_2$, we obtain
$$
-3\,\CAA'\cdot\CBB^2+\CBB^3=O(\delta_1).
$$

Note that $\delta_1^2\delta_2>g\delta_1$ by assumption. 
Then it remains to prove
\begin{equation*} 
\CLL_0^3=O(\delta_1^2\delta_2). 
\end{equation*}
The coefficients of $\delta_1^3$ and $\delta_1^2$ in $\CLL_0^3$ are $\ol\alpha_1^3$ and 
$3\, \ol\alpha_1^2\cdot (\ol\Delta-\ol\alpha_1-\ol\alpha_2)$ respectively. 
We claim that $\ol\alpha_1^3= \ol\alpha_1^2\cdot (\ol\Delta-\ol\alpha_1-\ol\alpha_2)=0$.
With the claim, the next largest monomial of $\delta_1,\delta_2$ is $\delta_1^2\delta_2$ since $\delta_1>\delta_2$ and $\delta_1\delta_2>g$. 
This proves the theorem.

It remains to prove the claim that $\ol\alpha_1^3= \ol\alpha_1^2\cdot (\ol\Delta-\ol\alpha_1-\ol\alpha_2)=0$. This is just a special case of the following basic result for the map $\tilde{p}_1:\CX\to\CCC$. 

\begin{lem} \label{general proj formula}
Let $f:\CX\to \CY$ be a morphism of arithmetic varieties over $\ZZ$.
Let $r = \dim \CY$ and $r + n = \dim\CX$.
Let $\eta\in \CY$ be the generic point and $\CX_\eta$ be the generic fiber of $f$. 
Let $\CLL_1, \dots, \CLL_n$ be hermitian line bundles on $\CX$ and $\CMM_1,\dots, \CMM_r$ be hermitian line bundles on $\CY$. 
Then  
$$
\CLL_1 \cdots \CLL_n \cdot f^*\CMM_1 \cdots f^* \CMM_r
=(\CL_{1,\eta} \cdots \CL_{n,\eta})  (\CMM_1 \cdots \CMM_r).
$$
Here the expressions in the three brackets are intersection numbers on $\CX,\CX_\eta,\CY$ respectively.
If $f$ is not surjective, then $\CX_\eta$ is empty and the intersection number on it is 0 by convention.
\end{lem}
\begin{proof}
We use induction on the relative dimension $n$.
The case $n=0$ follows from the projection formula. 
Assume that the result holds for $n-1$, and then consider the case $n$. Take a global section $s$ of $\CL_n$, and write 
$\div(s)=\sum_{i=1}^m a_i \CX_i$ for prime divisors $\CX_i$ on $\CX$. 
Denote by $f_i:\CX_i\to\CY$ the restriction of $f$ to $\CX_i$. 
We have the induction formula
\begin{multline*}
\CLL_1 \cdots \CLL_n \cdot f^*\CMM_1 \cdots f^* \CMM_r
= \sum_{i=1}^m a_i\, \CLL_1|_{\CX_i} \cdots \CLL_{n-1}|_{\CX_i} \cdot f_i^*\CMM_1 \cdots f_i^* \CMM_r\\
-\sum_{i=1}^m a_i\,\int_{\CX_i(\CC)} \log\|s\| c_1(\CLL_1) \cdots c_1(\CLL_n)  f^*c_1(\CMM_1) \cdots f^* c_1(\CMM_r).
\end{multline*}
The integral equals 0 since $c_1(\CMM_1) \cdots c_1(\CMM_r)=0$ on $\CY(\CC)$ by dimension reason.
It follows that
\begin{multline*}
\CLL_1 \cdots \CLL_n \cdot f^*\CMM_1 \cdots f^* \CMM_r
= \sum_{i=1}^m a_i\, \CLL_1|_{\CX_i} \cdots \CLL_{n-1}|_{\CX_i} \cdot f_i^*\CMM_1 \cdots f_i^* \CMM_r.
\end{multline*}
Let us first simplify every term of the right-hand side. 

If $\CX_i$ is vertical in the sense that the image of $\CX_i\to \Spec \ZZ$ is a closed point represented by a prime number $p$, then
\begin{equation*}
\CLL_1|_{\CX_i} \cdots \CLL_{n-1}|_{\CX_i} \cdot f_i^*\CMM_1 \cdots f_i^* \CMM_r=
\CL_1|_{\CX_i} \cdots \CL_{n-1}|_{\CX_i} \cdot f_i^*(\CM_1|_{\CY_p}) \cdots 
f_i^* (\CM_r|_{\CY_p}). 
\end{equation*}
is 0 as $r>\dim \CY_p$. 
If $\CX_i$ is horizontal  in that $\CX_i\to \Spec \ZZ$ is surjective, then it is an arithmetic variety, and the induction hypothesis gives 
\begin{equation*}
\CLL_1|_{\CX_i} \cdots \CLL_{n-1}|_{\CX_i} \cdot f_i^*\CMM_1 \cdots f_i^* \CMM_r=
(\CL_{1,\eta} \cdots \CL_{n-1,\eta}\cdot \CX_{i,\eta})  (\CMM_1 \cdots \CMM_r). 
\end{equation*}
Putting these into the formula, we get the result for $n$. 
\end{proof}

\section{Step 2: lower bound of the height}

Recall that Theorem \ref{small section} constructs a small section
$s$ of $\CLL=d\CLL_0$. Write 
$$
\CLL=(d_1-d) \tilde p_1^* \ol\alpha+(d_2-d) \tilde p_2^* \ol\alpha+ \CO(d\ol\Delta),
$$
where 
$$d_1=\delta_1d, \quad 
d_2=\delta_2d.$$
The goal here is to give a lower bound of $h_{\CLL}(P)$ using the section $s$. 

Denote by $\CL$ the underlying line bundle of $\CLL$, and by $L=\CL_K$ the generic fiber as before. 
Recall that $\CCC$ and $\CX$ are respectively integral models of $C$ and $X=C^2$ over $O_K$. 
Recall that $P=(P_1,P_2)$ is a $K$-point of $X$.
Denote by $\CP$ (resp. $\CP_1$, $\CP_2$) the Zariski closure of 
$P$ (resp. $P_1$, $P_2$) in $\CX$ (resp. $\CCC$, $\CCC$).

If $s(P)\neq 0$,  we would simply have 
$$
h_{\CLL}(P)=\wh\deg(\CLL|_\CP)
\geq - \sum_{\sigma:K\to \CC} \log\|s(P)\|_\sigma
\geq - \sum_{\sigma:K\to \CC} \log\|s\|_{\sigma,\sup}.
$$
This would imply Theorem \ref{vojta2}
by the bound in Theorem \ref{small section}. 
However, in general, we cannot avoid the case $s(P)=0$ where the above argument does not work. We will solve this problem by carefully analyzing of the effect of the index (which we will define later) of $s$ at $P$ on the height.

For $i=1,2$, choose a local coordinate $t_i$ of $C_i$ at $P_i$, i.e. $t_i$ is a generator of the maximal ideal of the local ring of $C_i$ at $P_i$. Let $\CO_X$ denotes the structure sheaf of $X$. For every local regular function $f$ at $P$ (i.e. $f\in\CO_{X,P}$), let $a_{i_1,i_2}$ be the coefficient of $t_1^{i_1}t_2^{i_2}$ in the power series expansion of $f$ in terms of variables $t_1, t_2$. Although $a_{i_1,i_2}$ may depends on the choice of $(s_0, t_1, t_2)$, the pairs $(i_1,i_2)$ such that $a_{i_1,i_2}\neq 0$ does not depend on the choice of $(s_0, t_1, t_2)$. Therefore, we can define \emph{the index of $f$ at $P$ of weight $(d_1, d_2)$} by 
\begin{equation}\label{eq:index of regular function}
\Ind_P(f, (d_1,d_2))
=\min\left\{ \frac{i_1}{d_1}+ \frac{i_2}{d_2}:
a_{i_1,i_2}\neq0
\right\}. 
\end{equation}

Now consider the line bundle $L$ and its section $s$. Choose a local section $s_0$ of $L$ at $P$ such that $s_0(P)\neq 0$. Then $s/s_0$ is a local regular function at $P$, and we define \emph{the index of $s$ at $P$ of weight $(d_1, d_2)$} by
$$
\Ind_P(s, (d_1,d_2))
=\Ind_P\left(\frac{s}{s_0}, (d_1,d_2)\right).
$$
Note that $\Ind_P(s, (d_1,d_2))$ is independent of the choice of $s_0$.

Assume that the index $\Ind_P(s, (d_1, d_2))$
is obtained at the monomial $t_1^{e_1}t_2^{e_2}$ of the power series expansion of $s$ at $P$. 
The following is the main result of this section, whose current form does not require the previous assumptions on $(\delta_1,\delta_2,d, P_1, P_2)$ and $s$. 

\begin{thm} \label{height inequality}
Let $s$ be a nonzero global section of $\CL$ on $\CX$. Then
\begin{multline*}
h_{\CLL}(P)
\geq - \sum_{\sigma:K\to \CC} \log\|s\|_{\sigma,\sup}
-\frac{2g-2}{g}[K:\QQ]\left(e_1 |P_1|^2
+e_2 |P_2|^2\right)\\
+O(d_1+d_2+d+e_1+e_2). 
\end{multline*}
Here the coefficients in the error terms are independent of $(\delta_1, \delta_2, d, P_1,P_2, s)$. 
\end{thm}

\subsection{The hermitian line bundle by restriction}

For $i=1,2$, denote by $\CPP_i=(\CP_i, g_{\CP_i})$ the arithmetic divisor on $\CX$ with underlying divisor $\CP_i$ and Arakelov's Green function $g_{\CP_i}$. This induces hermitian line bundles $\CO(\CPP_i)$ on $\CX$. 

Define a hermitian line bundle $\CMM$ on $\Spec O_K$ by 
$$
\CMM:=\big(\CLL-e_1\, \wt p_1^*\CO(\CPP_1)- e_2\, \wt p_2^* \CO(\CPP_2)\big)|_{\CP}.
$$
Here $\wt p_i:\CX\to \CCC$ is the morphism extending the projection $p_i:X\to C$.

Denote by $\ol\omega_{\CCC/O_K}$ the relative dualizing sheaf of $\CCC$ over $O_K$ endowed with the Arakelov metric $\VCV_\Ar$ as in \cite[\S4]{Ara74}. 
 Then we obtain the following result immediately. 

\begin{lem}
\begin{align*}
h_{\CLL}(P)
=&\, \wh\deg(\CMM)
-e_1 h_{\ol\omega_{\CCC/O_K}}(P_1)
-e_2 h_{\ol\omega_{\CCC/O_K}}(P_2)\\
=&\,\wh\deg(\CMM)
-\frac{2g-2}{g}[K:\QQ](e_1 |P_1|^2
+e_2 |P_2|^2)
+O(e_1+e_2).
\end{align*}
\end{lem}
\begin{proof}
By definition,
$$
\wh\deg(\CMM)
=\CLL\cdot \CP- e_1\, \wt p_1^*\CO(\CPP_1)\cdot \CP
- e_2\, \wt p_2^*\CO(\CPP_2)\cdot \CP.
$$
Note that
$$
\wt p_i^*\CO(\CPP_i)\cdot \CP
=\CO(\CPP_i)\cdot \CP_i
=-\ol\omega_{\CCC/O_K}\cdot\CP_i = -[K:\QQ]\frac{2g-2}{g}|P_i|^2 + O(1),
$$
where the first equality follows from the projection formula for $\wt p_i:\CX\to \CCC$, 
 the second equality follows from the arithmetic adjunction formula, 
 and the third equality is proved in Theorem \ref{mumford}. 
\end{proof}

\subsection{The section of the line bundle}

The underlying line bundle $\CM$ of $\CMM$ on $\Spec O_K$ is given by 
$$
\CM=\big(\CL-e_1\, \wt p_1^*\CO(\CP_1)- e_2\, \wt p_2^* \CO(\CP_2)\big)|_{\CP}.
$$
It can be identified with a locally free $O_K$-module of rank one.
We will define a nonzero section $s^\circ$ of $\CM$ from the section $s$ of $\CL$.

Start with the generic fiber 
$$
M=\big(L-e_1\,  p_1^*\CO(P_1)- e_2\,  p_2^* \CO(P_2)\big)|_{P}.
$$
As in the definition of the index, let $s_0$ be a local section of $L$ at $P$ with $s_0(P)\neq 0$, and for $i=1,2$ let $t_i$ be a local coordinate of $C_i$ at $P_i$. 
Then we have a power series expansion 
$$s=s_0(a_{e_1,e_2}t_1^{e_1}t_2^{e_2}+\dots).$$
Via the canonical injection $\CO(-P_i)\to \CO_C$, the local section $t_i$ of $\CO_C$ corresponds to a local section $t_i^\dagger$ of the line bundle $\CO(-P_i)$ with $t_i^\dagger(P_i)\neq 0$.
Finally, we obtain a nonzero element $s^\circ$ of $M$ by setting 
$$
s^\circ=(a_{e_1,e_2} s_0) \otimes (t_1^\dagger)^{\otimes e_1} \otimes (t_2^\dagger)^{\otimes e_2}
$$
according to the decomposition 
$$
M=\big(L+e_1\,  p_1^*\CO(-P_1)+ e_2\,  p_2^* \CO(-P_2)\big)|_{P}.
$$

\begin{lem}
The element $s^\circ\in M$ is independent of the choice of $(s_0,t_1,t_2)$, and lies in the lattice $\CM$ of $M$.
\end{lem}

\begin{proof}
We first prove the first statement. 
If we change $s_0$ to $s_0'=us_0$ for an invertible local regular function $u$ at $P$, then  the expansion
$$s=s_0(a_{e_1,e_2}t_1^{e_1}t_2^{e_2}+\dots)$$
becomes
$$s=s_0'(a_{e_1,e_2}'t_1^{e_1}t_2^{e_2}+\dots)$$
with $a_{e_1,e_2}'=u(0)^{-1}a_{e_1,e_2}$.
Thus $s^\circ$ does not change. 

Similarly, if we change $(t_1,t_2)$ to $(t_1',t_2')=(u_1t_1,u_2t_2)$ for invertible local regular functions $u_1, u_2$ at $P$, then the process does not change $s^\circ$ either.
This proves that $s^\circ\in M$ is independent of the choice of $(s_0,t_1,t_2)$.

To prove the second statement, we first review some geometric properties.
Since $\CCC$ is regular, the sections $\CP_1, \CP_2$ lie in the smooth locus $\CCC^\sm$ of $\CCC$ above $O_K$. 
See \cite[\S9.1, Cor. 1.32]{Liu02} for this well-known fact. 
By assumption, the natural map $\phi:\CX\to \CCC\times_{O_K}\CCC$ is an isomorphism above $\CCC^\sm\times_{O_K}\CCC^\sm$, so it is an isomorphism from a neighborhood of $\CP$ in $\CX$ to a neighborhood of $\CP_1\times_{O_K}\CP_2$ in 
$\CCC\times_{O_K}\CCC$. 

Now we prove the second statement that $s^\circ$ lies in $\CM$.
It suffices to prove that $s^\circ$ lies in every localization
$$\CM_{(\wp)}=\CM\otimes_{O_K} O_{K,(\wp)},$$ 
where $O_{K,(\wp)}$ denotes the local ring of $O_K$ at a prime ideal $\wp$. 
This will be done by an integral version of the definition of $s^\circ$. 
For $i=1,2$, let $\CP_{i,\wp}$ (resp. $\CP_\wp$) be a closed point of $\CP_i$ (resp. $\CP$) above $\wp$.  
Since $\CP_i$ lies in the smooth locus of $\CCC$, the maximal ideal of the local ring 
$\CO_{\CCC, \CP_{i,\wp}}$ can be generated by a single element, which we denote by $t_i$.
Moreover, the completion of $\CO_{\CCC, \CP_{i,\wp}}$ is isomorphic to the formal power series ring $O_{K,(\wp)}[[t_i]]$. 
Via the morphism $\phi$, the completion of $\CO_{\CX, \CP_{\wp}}$ is isomorphic to the formal power series ring $O_{K,(\wp)}[[t_1,t_2]]$. 
 
Let $s_0$ be a generator of the stalk of $\CL$ at $\CP_\wp$. 
Then $s/s_0$ lies in $\CO_{\CX, \CP_{\wp}}$, and thus has a power series expansion
$$
s/s_0=\sum_{i_1,i_2}a_{i_1,i_2} t_1^{i_1}t_2^{i_2},\quad  a_{i_1,i_2}\in O_{K,(\wp)}.
$$
In particular, we have $a_{e_1,e_2}\in O_{K,(\wp)}$.
It follows that 
$$
s^\circ=(a_{e_1,e_2} s_0) \otimes (t_1^\dagger)^{\otimes e_1} \otimes (t_2^\dagger)^{\otimes e_2}
$$
lies in 
$$
\CM=\big(\CL+e_1\,  \wt p_1^*\CO(-\CP_1)+ e_2\,  \wt p_2^* \CO(-\CP_2)\big)|_{\CP}.
$$
This completes the proof.
\end{proof}

\subsection{The degree of the line bundle}

In the above paragraphs, we have constructed a nonzero section 
$$
s^\circ=(a_{e_1,e_2} s_0) \otimes (t_1^\dagger)^{\otimes e_1} \otimes (t_2^\dagger)^{\otimes e_2}
$$
of the underlying line bundle $\CM$ of
$$
\CMM=\big(\CLL-e_1\, \wt p_1^*\CO(\CPP_1)- e_2\, \wt p_2^* \CO(\CPP_2)\big)|_{\CP}.
$$

By definition, for every embedding $\sigma:K\to\CC$, the metric
$$
\|s^\circ\|_\sigma
= |a_{e_1,e_2}|_\sigma\cdot  \|s_0(P)\|_\sigma \cdot \|t_1^\dagger(P_1)\|_\sigma^{e_1} \cdot \|t_2^\dagger(P_2)\|_\sigma^{e_2}.
$$
Here the metric 
$$
\|t_i^\dagger(P_i)\|_\sigma= 
\|t_i\cdot 1^\dagger(P_i)\|_\sigma= 
\lim_{Q\to P_i}(|t_i(Q)|_\sigma\cdot \|1^\dagger(Q)\|_\sigma)= 
\lim_{Q\to P_i}(|t_i(Q)|_\sigma\cdot e^{g_{\CP_i,\sigma}(Q)}),
$$
where $Q$ is a point on $C_\sigma(\CC)$ converging to $P_i$.
By the definition of the Arakelov metric $\|\cdot\|_\Ar$ on $\ol\omega_{\CCC/O_K}$, we have exactly
$$
\|t_i^\dagger(P_i)\|_\sigma= \|(\d t_i)(P_i)\|_{\Ar,\sigma}. 
$$
It follows that 
$$
\|s^\circ\|_\sigma
= |a_{e_1,e_2}|_\sigma\cdot  \|s_0(P)\|_\sigma \cdot \|(\d t_1)(P_1)\|_{\Ar,\sigma}^{e_1} \cdot \|(\d t_2)(P_2)\|_{\Ar,\sigma}^{e_2}.
$$

Now we are ready to estimate the degree of $\CMM$.
By definition, 
$$
\wh\deg(\CMM)=\log\# (\CM/sO_K)- \sum_{\sigma:K\to \CC} \log\|s^\circ\|_{\sigma}
\geq - \sum_{\sigma:K\to \CC} \log\|s^\circ\|_{\sigma}.
$$
It follows that 
\begin{multline}\label{eq:lower bound of degree}
\wh\deg(\CMM)
\geq - \sum_{\sigma:K\to \CC} 
\Big(
\log |a_{e_1,e_2}|_\sigma 
+\log \|s_0(P)\|_\sigma\\
+e_1\log \|(\d t_1)(P_1)\|_{\Ar,\sigma}
+e_2\log \|(\d t_2)(P_2)\|_{\Ar,\sigma}
\Big).
\end{multline}

Denote by  
$$D_{r_i}=\{Q\in C_\sigma(\CC): |t_i(Q)|\leq r_i\}$$
the closed disc of radius $r_i$ in $C_\sigma(\CC)$ with center $P_i$. 
We can choose $r_1,r_2$ small enough so that $s_0$ is holomorphic and everywhere non-vanishing at a neighborhood of the poly-disc 
$$D=D_{r_1,r_2}=D_{r_1}\times D_{r_2}$$ 
in $X_\sigma(\CC)$.  
Then $f=s/s_0$ is a holomorphic function on $D$.
By Cauchy's integration formula, 
$$
a_{e_1,e_2} = \frac{1}{(2\pi \i)^2}\int_{\partial D_{r_1}} \int_{\partial D_{r_2}} \frac{f(z_1,z_2)}{z_1^{e_1+1}z_2^{e_2+1}}\d z_2 \d z_1.
$$
It follows that 
$$
|a_{e_1,e_2}|_\sigma \leq r_1^{-e_1}r_2^{-e_2} |f|_{D,\sup}.
$$
Here $|f|_{D,\sup}$ denotes the maximal value of $|f|$ on $D$, and we will use other similar notation.
For example, we have 
$$
|f|_{D,\sup}=\|s/s_0\|_{D,\sup}
\leq \|s\|_{D,\sup}/\|s_0\|_{D,\inf}
\leq \|s\|_{\sigma, \sup}/\|s_0\|_{D,\inf}.
$$

Therefore, \eqref{eq:lower bound of degree} becomes
\begin{multline*}
\wh\deg(\CMM)
 \geq - \sum_{\sigma:K\to \CC} 
\Big(
\log \|s\|_{\sigma, \sup} 
+\log \frac{\|s_0(P)\|_\sigma}{\|s_0\|_{D,\inf}} 
-e_1\log r_1- e_2\log r_2 \\
+e_1\log \|(\d t_1)(P_1)\|_{\Ar,\sigma}
+e_2\log \|(\d t_2)(P_2)\|_{\Ar,\sigma}
\Big).
\end{multline*}

\subsection{The uniform choice}

Now we are going to choose $(s_0, (t_1,t_2))$ carefully to get a uniform bound with respect to $P$.
Consider the set of quadriples $(U_1\times U_2, V, (z_1,z_2), s_0)$, where 
\begin{enumerate}[(1)]
\item $U_1, U_2$ are open subsets of $C_\sigma(\CC)$ under the complex topology;
\item $z_1\colon U_1\to \{w\in \CC:|w|<3\}$ and $z_2\colon U_2\to \{w\in \CC:|w|<3\}$ 
are biholomorphic maps;
\item $V=\{Q\in U_1\times U_2: |z_1(Q)|< 1, \ |z_2(Q)|< 1\}$ is a poly-disc in $U$; 
\item $s_0$ is a holomorphic section of $L$ on $U$ which is everywhere non-vanishing on $U$. 
\end{enumerate}
By compactness, we can choose finitely many quadriples $(U_1\times U_2, V, (z_1,z_2),\linebreak s_0)$ such that these poly-discs $V$ cover $X_\sigma(\CC)$. For every point $P=(P_1,P_2)\in X_\sigma(\CC)$, there exists some chosen quadriple $(U_1\times U_2, V, (z_1,z_2), s_0)$ such that $P\in V$.
Then we set $(t_1,t_2)= (z_1-z_1(P_1),z_2-z_2(P_2))$. 
Take the poly-disc 
$$
D=\{Q\in U: |t_1(Q)|\leq 1,\ |t_2(Q)|\leq 1\},
$$
which is contained in 
$$
\{Q\in U_1\times U_2: |z_1(Q)|< 2, \ |z_2(Q)|< 2\}.
$$
By the currently chosen $(s_0,(t_1, t_2))$ for $P$, we have the bound
\begin{align*}
\wh\deg(\CMM)
 \geq &  - \sum_{\sigma:K\to \CC} 
\Big(
\log \|s\|_{\sigma, \sup} 
+\log \frac{\|s_0(P)\|_\sigma}{\|s_0\|_{D,\inf}} 
-e_1\log r_1- e_2\log r_2 \\
& +e_1\log \|(\d t_1)(P_1)\|_{\Ar,\sigma}
+e_2\log \|(\d t_2)(P_2)\|_{\Ar,\sigma}
\Big).
\end{align*}
We see that $r_1=r_2=1$, and
$$\log\|s_0(P)\|_\sigma,\quad  \log\|s_0\|_{D,\inf},\quad 
\log \|(\d t_1)(P_1)\|_{\Ar,\sigma},\quad 
\log\|(\d t_2)(P_2)\|_{\Ar,\sigma}$$
are all bounded from both above and below uniformly with respect to $P$. 
Moreover, as $(d_1,d_2,d)$ varies,
we can choose the section $s_0$ of 
$$
L=(d_1-d)  p_1^* \alpha+(d_2-d)  p_2^* \alpha+ \CO(d\Delta)
$$
on every $U$ as tensor powers of fixed sections of the individual line bundles 
$p_1^* \alpha, p_2^* \alpha, \CO(d\Delta)$.
This makes 
$\log\|s_0(P)\|_\sigma$
and 
$\log\|s_0\|_{D,\inf}$
grows as $O(d_1+d_2+d)$.
As a consequence, we have 
\begin{align*}
\wh\deg(\CMM)
 \geq  - \sum_{\sigma:K\to \CC} 
\log \|s\|_{\sigma, \sup} 
+O(d_1+d_2+d)+O(e_1+e_2).
\end{align*}

\section{Step 3: upper bound of the index}

The goal of this section is to provide an upper bound of the index used in Theorem \ref{height inequality}. This lies in the framework of the classical Dyson's lemma.
Our exposition follows from \cite{Voj89a}, which in turn is inspired by \cite{EV84, Vio85}.
Let us first introduce the framework of this section, which is slightly different from that of the previous sections.

Let $K$ be an algebraically closed field of characteristic $0$. 
Let $C_1, C_2$ be smooth projective curves of genera $g_1,g_2$ over $K$ respectively. 
Let $X=C_1\times C_2$ be the product, and $p_1\colon X\to C_1,$ $p_2\colon X\to C_2$ be the projection morphisms.
Denote by $F_1$ (resp. $F_2$) a closed fiber of $p_1$ (resp. $p_2$), which will be used for computing intersection numbers.
Let $D$ be a nonzero effective divisor on $X$. 
Denote by $e(D)$ the maximum of $\ord_Y D$ over all prime divisors $Y$ of $X$ which are not fibers of $p_1$ or $p_2$.

Let $b_1,b_2$ be positive integers. For every point $P=(P_1,P_2)\in X(K)$, we define the 
\emph{index of $D$ at $P$ of weight $(b_1, b_2)$} to be 
$$
\Ind_P(D, (b_1, b_2))
:= \Ind_P(f, (b_1, b_2)),
$$
where $f$ is a local generator of the ideal sheaf $\CI_D$ at $P$, i.e. the stalk $\CI_{D,P}$, and the index on the right is defined by \eqref{eq:index of regular function} (here we view $f$ as a local regular function in $\CO_{X,P}$ via the inclusion $\CI_D\subset \CO_X$). Note that $\Ind_P(D,(b_1, b_2))$ does not depend on the choice of the local generator $f$. 

Let $Q_1,\dots, Q_m$ be distinct $K$-points on $X=C_1\times C_2$.
For $i=1,2$, denote by 
$$m_i=m_i(Q_1,\dots, Q_m)$$ 
the number of different points in $p_i(Q_1), \dots, p_i(Q_m)$. For $k=1,\dots, m$, let
$$I_k:=\Ind_{Q_k}(D, (b_1, b_2))$$
denote the index.

Finally, denote a function $V:[0,\infty)\to \RR$ by
$$
V(a):=\vol\Big(\{(x,y)\in\RR^2: 0\leq x\leq 1, \ 0\leq y\leq 1, \ x+y\leq a\}\Big).
$$
Here the right-hand side is the area in the Euclidean plane.

Finally, the main theorem of this section is the following theorem from \cite[Thm. 0.4]{Voj89a}. 

\begin{thm}[generalized Dyson's lemma]\label{dyson}
Assume that $m_1=m_2=m$. 
Assume that $b_1\leq D\cdot F_2$ and $b_2\leq D\cdot F_1$ are positive real numbers.
Then
$$
\sum_{k=1}^m V(I_k) \leq \frac{D\cdot D}{2b_1b_2}
+e(D)\frac{D\cdot F_1}{2b_1b_2} \max\{2g_1-2+m_1,0\}. 
$$
\end{thm}

The intuition for the inequality is as follows. Write $D=\div(s)$ for a nonzero section $s$ of a line bundle $L$ on $X$. Fix $X, L, Q_1, \dots, Q_m, b_1, b_2$, but vary $s$ in $H^0(X, L)$. Fix positive real numbers $I_1',\dots, I_m'$ in the interval $(0,1)$.
Consider conditions for the existence of a section $s$ satisfying $\Ind_{Q_k}(\div(s), (b_1, b_2))\geq I_k'$. These conditions become linear equations on the section $s$ at $Q_1, \dots, Q_m$, and the number of linear equations is roughly $\sum_{k=1}^m b_1b_2 V(I_k')$. If these linear equations are linearly independent, then the existence of $s$ is equivalent to 
$$\sum_{k=1}^m b_1b_2 V(I_k')< h^0(X, L)+\text{error}.$$
If $L$ is sufficiently ample, we expect 
$$h^0(X, L)=\ds\frac12 L^2+\text{error}. $$
The combination of these two relations explains the inequality in Dyson's lemma.

\subsection{First reduction: remove fibers}

We first reduce Theorem \ref{dyson} to the following case where $D$ contains no fibers of $p_1$ or $p_2$ and the assumptions $m_1=m_2=m$, $b_1\leq D\cdot F_2$, $b_2\leq D\cdot F_1$ are removed.

\begin{thm} \label{dyson2}
Assume that $D$ does not contain any fiber of $p_1$ or $p_2$.
Then
$$
\sum_{k=1}^m \frac12 I_k^2 \leq \frac{D\cdot D}{2b_1b_2}
+e(D)\frac{D\cdot F_1}{2b_1b_2} \max\{2g_1-2+m_1,0\}. 
$$
\end{thm}

Let us prove Theorem \ref{dyson} by Theorem \ref{dyson2}.
Let $D$ be as in Theorem \ref{dyson}. 
If $D$ contains a fiber of $p_1$ or $p_2$ which does not pass through any of $Q_1,\dots, Q_m$, then removing this fiber from $D$ does not change the left-hand side of the inequality of the theorem while decrease the right-hand side of the inequality. 
Therefore, we can assume that $D$ does not contain such fibers.

Write 
\begin{multline*}
 D=D'+\sum_{k=1}^{m} \left(x_k\left(p_1^{-1}\left(p_1(Q_k)\right)\right) +y_k\left(p_2^{-1}\left(p_2(Q_k)\right)\right)\right)\\
 = D'+\sum_{k=1}^{m} \left(x_k p_1(Q_k)\times C_2 +y_kC_1\times p_2(Q_k)\right),
\end{multline*}
where the effective divisor $D'$ on $X$ does not contain any fiber of $p_1$ or $p_2$,
and the multiplicities $x_1,\dots, x_m, y_1, \dots, y_m$ are 
non-negative integers.
Since $m_1 = m_2 = m$, i.e. $p_1(Q_1),\ldots,p_1(Q_m)$ are distinct and $p_2(Q_1),\ldots,p_2(Q_m)$ are also distinct, it is easy to see that the indices $I_k=\Ind_{Q_k}(D, (b_1,b_2))$ and $I_k'=\Ind_{Q_k}(D', (b_1,b_2))$ satisfy
$$
I_k=I_k'+ \frac{x_k}{b_1} +\frac{y_k}{b_2}.
$$

By Theorem \ref{dyson2}, 
$$
\sum_{k=1}^m \frac12 I_k'^2 \leq \frac{D'\cdot D'}{2b_1b_2}
+e(D')\frac{D'\cdot F_1}{2b_1b_2} \max\{2g_1-2+m_1,0\}. 
$$
Note $e(D)=e(D')$, and $V(a)\leq a^2/2$ for all $a\geq0$. 
By taking difference, it suffices to prove  
\begin{multline*}
 \sum_{k=1}^m \left(V(I_k)-V(I_k')\right) \leq \frac{D\cdot D-D'\cdot D'}{2b_1b_2}\\
 +e(D)\frac{(D-D')\cdot F_1}{2b_1b_2} \max\{2g_1-2+m_1,0\}.
\end{multline*}
By the expression of $D$ in terms of $D'$, we can view both sides of the inequality as functions in non-negative variables $x_1,\dots, x_m, y_1, \dots, y_m$. 
When the variables are 0, the inequality is an equality. 
Therefore, it suffices to prove the corresponding inequality for the partial derivatives of both sides.
Note that the derivative $V'(a)$ exists and 
$V'(a)\leq 1$ for all $a\geq 0$.
Then 
$$
\frac{\partial}{\partial x_k}
\sum_{i=1}^m \left(V(I_i)-V(I_i')\right)
=\frac{\textnormal{d}}{\dif x_k} V(I_k)
=\frac{1}{b_1} V'(I_k)
\leq \frac{1}{b_1}.
$$
On the other hand, 
$$
\frac{\partial}{\partial x_k}
\left(\frac{D\cdot D-D'\cdot D'}{2b_1b_2}\right)
=\frac{2\, D\cdot F_1}{2b_1b_2}
\geq \frac{1}{b_1}.
$$
Here we have used the assumption $D\cdot F_1\geq b_2$.
This proves the inequality for the partial derivative with respect to $x_k$, and that for $y_k$ is similar. 
This completes the first reduction process.

\subsection{Second reduction: trivialize weight}

Now we reduce Theorem \ref{dyson2} to the case $b_1=b_2=1$ by the covering trick of \cite[\S3]{Voj89a}. The proof of \cite{Voj89a} applies the Kodaira--Parshin construction to construct ramified covering of $C_i$. Here we use a topological method instead, which works best over $\BC$.

\begin{lem} \label{covering}
Let $C$ be a smooth complex projective curve of genus $g$. 
Let $n>1$ be a positive integer.
Let $S$ be a finite set of points of $C$.
In the case $g=0$, we further assume that $|S|\geq2$.
Then there is a smooth projective complex curve $C'$ with  
a finite and flat morphism $f:C'\to C$ such that any point $P'\in C'$ is ramified over 
$C$ if and only if $f(P')\in S$, and in that case, the ramification index of $P'$ is exactly $n$.
\end{lem}

\begin{proof}

If $S=\emptyset$, then $g>0$, and the result is elementary in algebraic topology. Roughly speaking, an unramified cyclic covering of $C$ of order $n$ is given by a surjection $\pi_1(C)\to \ZZ/n\ZZ$, or equivalently a surjection $H_1(C,\ZZ)\to \ZZ/n\ZZ$ by taking abelianization.
There are plenty of such maps as $H_1(C,\ZZ)\cong \ZZ ^{2g}$. 

Now we assume that $|S|\geq 1$. By taking an unramified covering $h\colon C_2\to C$ of degree 2, and replacing $(C, S)$ by $(C_2,h^{-1}(S))$, we can assume further that $|S|\geq 2$. 
By factoring $n$, we can assume that $n$ is a prime number. 

For every $P\in S$, take a closed disc $D_P$ in $C$ with center $P$ of positive radius such that these discs are disjoint. There is a basic exact sequence
$$
0\lra \ZZ \stackrel{\iota_1}{\lra} \bigoplus_{P\in S} H_1(D_P\setminus P, \ZZ) \stackrel{\iota_2}{\lra} H_1(C\setminus S, \ZZ)
\lra H_1(C, \ZZ) \lra 0.
$$
Here $H_1(D_P\setminus P, \ZZ)\cong \ZZ$ is generated by the boundary of $D_P$, and the arrow $\iota_1$ is the diagonal map under this identification. 
Then $H_1(C\setminus S, \ZZ)\cong \ZZ ^{2g+|S|-1}$.

Let $H_1(C\setminus S, \ZZ)\to \ZZ/n\ZZ$ be a surjection whose kernel does not contain the image of $H_1(D_P\setminus P, \ZZ)$ for any $P\in S$. 
The surjection induces an unramified cyclic covering $f^\circ:U'\to C\setminus S$ of degree $n$. 
Let $C'$ be the unique smooth compactification of $U'$. Then we obtain a possibly ramified covering $f:C'\to C$. 
We will prove that this covering satisfies the requirement. 

By the choice, the surjection $H_1(C\setminus S, \ZZ)\to \ZZ/n\ZZ$ does not factor through $H_1(C\setminus (S\setminus P), \ZZ)\to \ZZ/n\ZZ$ for any $P\in S$. 
In other words, the map $f:C'\to C$ is ramified above $P$. As $f^\circ$ is Galois, the Galois action extends to $C'$. Then the ramification indices of all points in $f^{-1}(P)$ are equal and divides $n$, so they must be equal to $n$ by the assumption that $n$ is prime. 
This finishes the proof.
\end{proof}

Now we are ready to reduce Theorem \ref{dyson2} to the case $b_1=b_2=1$. 
It suffices to reduce it to the case $b_1=b_2$, since in the later case both sides of the inequality are quadratic in $b_1^{-1}$. 

If $g_i=0$ and $m_i=1$ for some $i=1,2$, then we can just add a point $Q_{m+1}$ of $X\setminus |D|$ to the set $\{Q_1,\dots, Q_m\}$ such that $m_1\geq 2$ and $m_2\geq 2$. This does not change both sides of the inequality.
Therefore, we can assume that for every $i=1,2$, we have either $g_i>0$ or $m_i\geq2$.

By the Lefschetz principle, we can assume that $(C, D, Q_1,\dots, Q_m)$ are defined over a finitely generated fields of characteristic 0, and by base change we can assume that they are actually defined over $\CC$. This reduce the problem to the case $K=\CC$.

Now we are ready to apply Lemma \ref{covering}. 
Set $(n_1,n_2)=(b_2,b_1)$. 
By the lemma, for $i=1,2$, there is a covering $f_i:C_i'\to C_i$ which is unramified outside $\{p_i(Q_1),\dots, p_i(Q_m)\}$ and the ramification index of points above $\{p_i(Q_1),\dots, p_i(Q_m)\}$ are equal to $n_i$.
Let $X=C_1'\times C_2'$, and let $f\colon X'\to X$ the induced map. 
Denote by $Q_1',\dots, Q_{m'}'$ the points of $f^{-1}(Q_1,\dots, Q_m)$.

Now we claim that Theorem \ref{dyson2} for $(X,(Q_1',\dots, Q_{m'}'), f^*D, (b_1n_1,\linebreak b_2n_2))$
implies that for $(X,(Q_1,\dots, Q_{m}), D, (b_1,b_2))$.
In fact, we have the relation 
$$
\Ind_{Q'}(D', (b_1n_1,b_2n_2))=\Ind_{Q_k} (D, (b_1,b_2))
$$
for any $Q'\in f^{-1}(Q_k)$. 
The Riemann--Hurwitz formula gives the relation
$$
2g_1'-2+m_1'= \deg(f_1)(2g_1-2+m_1),
$$
where $g_1'$ denotes the genus of $C_1'$,  $m_1'$ denotes the number of distinct points in $p_1'(Q_1',\dots, Q_{m'}')$, and $p_1':X'\to C_1'$ denotes the projection. 
By these, we can check that the left-hand side (resp. right-hand side) of the inequality in Theorem \ref{dyson2} for $(X,(Q_1',\dots, Q_{m'}'), f^*D, (b_1n_1,b_2n_2))$ is equal to 
$\deg(f_1)\deg(f_2)/(n_1,n_2)$ times that for $(X,(Q_1,\dots, Q_{m}), D, (b_1,b_2))$. 
This finishes the reduction by $b_1n_1=b_2n_2$.

\subsection{Proof of the essential case}

In the above paragraphs, we reduced Theorem \ref{dyson2} to the case $b_1=b_2=1$. 
Now we prove this case. 

Denote by $\psi:X'\to X$ the blowing-up of $X$ along $\{Q_1,\dots, Q_m\}$. 
Denote by $E_1,\dots, E_m$ the exceptional divisors above $Q_1,\dots, Q_m$ respectively.
Write 
$$
D=a_1D_1+\dots+a_r D_r 
$$
in terms of distinct prime divisors $D_1,\dots, D_r$ of $X$.
Here the coefficients $a_1,\dots, a_r$ are strictly positive.
For $j=1,\dots, r$, denote by $D_j'$ the strict transform of $D_j$. 
Then the strict transform of $D$ is just 
$$
D'=a_1D_1'+\dots+a_r D_r'.
$$

By a direct calculation, the index $I_k=\Ind_{Q_k}(D, (1,1))$ is the multiplicity of $E_k$ in the pull-back $\psi^*D$. 
In other words, we have 
$$
\psi^*D=D'+I_1 E_1+\dots+I_m E_m.
$$
By the projection formula for $\psi$, the intersection number $\psi^*D\cdot E_k=0$. 
It follows that 
$$
D'^2=(\psi^*D-I_1 E_1+\dots+I_m E_m)^2=D^2-I_1^2-\dots-I_m^2.
$$
Then 
$$
I_1^2+\dots+I_m^2=D^2-D'^2.
$$
Therefore, Theorem \ref{dyson2} (for $b_1=b_2=1$) is equivalent to 
$$
D'^2 
+e(D)(D\cdot F_1) \max\{2g_1-2+m_1,0\}
\geq 0. 
$$
This is further equivalent to 
$$
D'^2 
+e(D')(D'\cdot \psi^*F_1) \max\{2g_1-2+m_1,0\}
\geq 0. 
$$
Here $e(D')=e(D)=\max\{a_1,\dots, a_r\}$.

We first prove the inequality holds for a single irreducible component. Namely, we have the following stronger statement, which is the reason for all other generalization.

\begin{thm} \label{essential dyson}
Let $Y$ be a integral closed curve of $X'$ which is not contained in fibers of $p_1':X'\to C_1$ or $p_2':X'\to C_2$. Then 
\begin{equation}\label{eq:essential dyson}
 Y^2 
 +(Y\cdot \psi^*F_1) (2g_1-2+m_1)
 \geq 0.
\end{equation}
\end{thm}

\begin{proof}
As in \cite{Voj89a}, the key of the proof is the adjunction formula and the Riemann--Hurwitz formula. 

Recall the adjunction formula
\begin{equation}\label{eq:adjunction formula}
 2h^1(\CO_Y)-2 = Y^2+Y\cdot \omega_{X'},
\end{equation}
which is well-known if $Y$ is smooth. The singular case is a combination of the expression
$$
\chi(\CO_Y)=\chi(\CO_{X'})-\chi(\CO_{X'}(-Y)),
$$
which comes from the additivity of the Euler characteristic $\chi$, 
and the Riemann--Roch theorem for the line bundle $\CO_{X'}(-Y)$ on $X'$. 

Let $f\colon \wt Y\to Y$ be the normalization of $Y$. Then we have the following exact sequence
$$
0\lra \CO_{Y} \lra f_* \CO_{\wt Y} \lra \CQ\lra 0.
$$
Since the quotient sheaf $\CQ=f_* \CO_{\wt Y}/\CO_{Y}$ is supported on the singular locus of $Y$, $H^1(Y, \CQ)=0$. 
Thus we see that
\begin{equation}\label{eq:h^1 of nomalization}
 h^1(\CO_Y) \geq h^1(\CO_{\wt Y}).
\end{equation}
By the Riemann--Hurwitz formula,
\begin{equation}\label{eq:Riemann-Hurwitz for normalization}
 2 h^1(\CO_{\wt Y})-2 \geq \deg(\wt Y/C_2) (2g_2-2)
 =(Y\cdot \psi^*F_2) (2g_2-2)
 =Y\cdot p_2'^*\omega_{C_2}.
\end{equation}
Then \eqref{eq:adjunction formula}, \eqref{eq:h^1 of nomalization}, and \eqref{eq:Riemann-Hurwitz for normalization} imply that
\begin{equation}\label{eq:adjunction inequality}
 Y^2+Y\cdot \omega_{X'} \geq Y\cdot p_2'^*\omega_{C_2}.
\end{equation}
Moreover, the canonical sheaf of blowing-up is given by
$$
\omega_{X'}=\psi^* \omega_{X}+\CO(E_1+\dots+E_m)
=p_1'^* \omega_{C_1}+p_2'^* \omega_{C_2}+\CO(E_1+\dots+E_m). 
$$
Thus \eqref{eq:adjunction inequality} becomes
$$
Y^2+ Y\cdot p_1'^* \omega_{C_1}+Y\cdot E_1+\dots+Y\cdot E_m\geq 0.
$$
In order to prove \eqref{eq:essential dyson}, it remains to prove that
$$
Y\cdot E_1+\dots+Y\cdot E_m \leq (Y\cdot \psi^*F_1)m_1 .
$$

Finally,
denote by $\{P_1,\dots, P_{m_1}\}$ 
the distinct points of $\{p_1(Q_1),\dots, p_1(Q_m)\}$. 
Then the divisor $F=p_1^{-1}(P_1)+\dots + p_1^{-1}(P_{m_1})$ passing through all points of
 $\{Q_1,\dots, Q_m\}$. 
By the blowing-up process, 
$\psi^*F-E_1-\dots-E_m$ is an effective divisor intersecting $Y$ properly. 
This gives 
$$Y\cdot (\psi^*F-E_1-\dots-E_m)\geq 0.$$
This completes the proof by the fact that $F$ is numerically equivalent to $m_1F_1$.
\end{proof}

Now we are ready to complete the proof of Theorem \ref{dyson2}. 
Recall that right before Theorem \ref{essential dyson}, we are left to prove that
$$
D'^2 
+e(D')(D'\cdot \psi^*F_1) \max\{2g_1-2+m_1,0\}
\geq 0. 
$$
Here 
$$
D'=a_1D_1'+\dots+a_r D_r', \quad  e(D')=\max\{a_1,\dots, a_r\}.
$$
As the proper intersection $D_j'\cdot D_{j'}'\geq 0$ for $j\neq j'$, we have
$$
D'^2\geq \sum_{j=1}^r a_j^2D_j'^2.
$$
Since $F_1$ is nef, we have
$$
e(D')(D'\cdot \psi^*F_1) \geq \sum_{j=1}^r a_j^2(D_j'\cdot \psi^*F_1). 
$$
Thus it suffices to prove that
\begin{equation}\label{eq:last inequality}
 \sum_{j=1}^r a_j^2D_j'^2
 + \sum_{j=1}^r a_j^2(D_j'\cdot \psi^*F_1) \max\{2g_1-2+m_1,0\}
 \geq 0.
\end{equation}
By Theorem \ref{essential dyson}, 
$$
D_j'^2+(D'_j\cdot \psi^*F_1) (2g_1-2+m_1)
\geq 0. 
$$
Sum over $j$, and then \eqref{eq:last inequality} is proved.

\begin{remark}
The proof only works for base fields of characteristic 0 due to the application of the Riemann--Hurwitz formula. In fact, Dyson's lemma does not hold for positive characteristics by the example in \cite[\S4]{Voj89a}.
\end{remark}

\section{Completion of the proof}

Now we are prepared to prove Theorem \ref{vojta3}. Assume the setting of this theorem.
Recall that
$$
\CLL=d\CLL_0=(d_1-d) \tilde p_1^* \ol\alpha+(d_2-d) \tilde p_2^* \ol\alpha+ \CO(d\ol\Delta),
$$
where 
$$d_1=\delta_1d, \quad 
d_2=\delta_2d.$$
By Theorem \ref{small section}, 
for sufficiently large and sufficiently divisible integer $d$,
there exists a nonzero element
$s\in \Gamma(\CX, \CL)$ satisfying 
$$
-\log \|s\|_{\sigma,\sup}
\geq O\left(\frac{d\delta_1^2\delta_2}{\delta_1\delta_2-g}\right)
=O\left(\frac{d\delta_1}{\delta_1\delta_2-g}\right)
$$
for every embedding $\sigma\colon K\to \CC$. 
Here we assume that $\Ind_P(s, (d_1,d_2))$ is attained at exponent $(e_1,e_2)$, and we use $g<\delta_1\delta_2<g+1/4$ to simplify the error term. 

Taking this section $s$ in Theorem \ref{height inequality}, we have
$$
h_{\CLL}(P)_\QQ
\geq 
-\frac{2g-2}{g}(e_1 |P_1|^2
+e_2 |P_2|^2)
+O\left(\frac{d\delta_1}{\delta_1\delta_2-g}\right)
+O(e_1+e_2).
$$
We need the following bound of the index terms. 

\begin{lem}\label{bound index}
Let $\delta_1,\delta_2$ be positive rational numbers satisfying
\begin{equation*} 
\delta_1>2g\delta_2,  \quad  g<\delta_1\delta_2<g+\frac14.
\end{equation*}
Then the exponent $(e_1,e_2)$ satisfies
$$
e_1+e_2<  c(\delta_1,\delta_2)\delta_1 d,
$$
and
$$
e_1|P_1|^2+e_2|P_2|^2
\leq  c(\delta_1,\delta_2)\cdot   (\delta_1|P_1|^2+\delta_2|P_2|^2)\cdot d. 
$$
Here 
$$
c(\delta_1,\delta_2)=\sqrt{\frac{2}{g}(\delta_1\delta_2-g)+(2g-1)\frac{\delta_2}{\delta_1}}<1.
$$
\end{lem}

By the lemma, the bound of the height becomes
$$
h_{\CLL}(P)_\QQ
\geq 
-\frac{2g-2}{ g} c(\delta_1,\delta_2)\cdot   (\delta_1|P_1|^2+\delta_2|P_2|^2)\cdot d
+ O\left(\frac{d\delta_1}{\delta_1\delta_2-g}\right).
$$
It follows that
$$
h_{\CLL_0}(P)_\QQ
\geq 
-\frac{2g-2}{ g} c(\delta_1,\delta_2)\cdot   (\delta_1|P_1|^2+\delta_2|P_2|^2)
+ O\left(\frac{\delta_1}{\delta_1\delta_2-g}\right),
$$
which proves Theorem \ref{vojta3}. 

Finally, we deduce Lemma \ref{bound index} from Theorem \ref{dyson} as follows.

\begin{proof}[Proof of Lemma \ref{bound index}]
Set $C_1=C_2=C$ and $(b_1,b_2)=(d_1,d_2)$ to compute the index 
$$I=\Ind_{P}(s,(d_1,d_2))=\Ind_{P}(\div(s),(d_1,d_2)).$$
Then Theorem \ref{dyson} gives
$$
V(I) \leq \frac{d_1d_2-gd^2}{d_1d_2}
+e(\div(s))\frac{1}{2d_1} (2g-1). 
$$
Here $e(\div(s))$ is the maximal multiplicity of irreducible components $Y$ of $D$ which are not fibers of $p_1$ or $p_2$. 
It follows that 
$$
d_2=L\cdot F_1 \geq e(\div(s)) D\cdot F_1 \geq e(\div(s)). 
$$
As a consequence, we have
$$
V(I) \leq \frac{d_1d_2-gd^2}{d_1d_2}
+\frac{(2g-1)d_2}{2d_1} 
=\frac{\delta_1\delta_2-g}{\delta_1\delta_2}
+\frac{(2g-1)\delta_2}{2\delta_1}. 
$$
Since $\delta_1>2g\delta_2$ and
$g<\delta_1\delta_2<g+1/4$,
$$
\frac{\delta_1\delta_2-g}{\delta_1\delta_2}
+\frac{(2g-1)\delta_2}{2\delta_1}
<\frac{1}{4g}
+\frac{2g-1}{2\cdot 2g}
=\frac12. 
$$

It follows that $V(I)< 1/2$, and thus 
$$V(I)=\frac{1}{2}I^2=\frac12\left(\frac{e_1}{d_1}+\frac{e_2}{d_2}\right)^2.
$$ 
As a consequence,
$$
\left(\frac{e_1}{d_1}+\frac{e_2}{d_2}\right)^2\leq 
\frac{2(\delta_1\delta_2-g)}{\delta_1\delta_2}+(2g-1)\frac{\delta_2}{\delta_1}. 
$$

By $\delta_1\delta_2>g$, we have 
$$
\frac{e_1}{d_1}+\frac{e_2}{d_2}< c(\delta_1,\delta_2), 
$$
and thus 
$$
e_1 \delta_2+e_2 \delta_1< c(\delta_1,\delta_2)\delta_1\delta_2 d. 
$$
As $\delta_1>\delta_2$, this implies 
$$
e_1+e_2<  c(\delta_1,\delta_2)\delta_1 d. 
$$
On the other hand, the trivial bound
$$
e_1|P_1|^2+e_2|P_2|^2
\leq (e_1 \delta_2+e_2 \delta_1)\cdot  (\delta_1|P_1|^2+\delta_2|P_2|^2) \cdot (\delta_1\delta_2)^{-1}
$$
gives 
$$
e_1|P_1|^2+e_2|P_2|^2
\leq  c(\delta_1,\delta_2)  d\cdot (\delta_1|P_1|^2+\delta_2|P_2|^2). 
$$
\end{proof}

\

{\footnotesize

\noindent Xinyi Yuan

\noindent Address: \emph{BICMR, Peking University, Haidian District, Beijing 100871, China}

\noindent Email: \emph{yxy@bicmr.pku.edu.cn}
}

\end{document}